\def \dd#1{{\bf#1}}

\def\cl#1{{\cal#1}}



\def\ouv#1{\smash{\mathop{#1}\limits^{\lower 1pt\hbox
{$\scriptscriptstyle\circ$}}}}

\def\hfl#1#2{\smash{\mathop{\hbox to 12mm{\rightarrowfill}}
\limits^{\scriptstyle#1}_{\scriptstyle#2}}}


\def\stit#1{\vskip 3mm plus 1mm minus 2mm {\bf{#1}}
		\smallskip}

\font\tir=cmbx10 at 12pt

\def\ref#1#2#3#4{{\bf #1}{\ #2}{\it ,\ #3}{,\ #4}\medskip}


\def \picture #1 by #2 (#3){\midinsert \centerline 
{\vbox to #2{\hrule width #1 heigth 0pt 
depth 0pt \null \vfill \special {picture #3}}}\endinsert }

\def\scaledpicture #1 by #2 (#3 scaled #4) {{
\dimen0 =#1 \dimen1 =$2
\divide \dimen0 by 1000 \multiply \dimen0 by #4
\divide \dimen1 by 1000 \multiply \dimen1 by #4
\picture \dimen0 by \dimen1 (#3 scaled $4)}}

\def\figure #1 #2 #3 {\midinsert \vglue 3mm 
{\vbox to #3 {\hrule width 6cm height 0cm depth 0cm \vfill
{\special {picture #1 scaled #2}}}}\vglue 2mm \endinsert}

\magnification=1200

\input psfig.sty

\overfullrule=0pt

{\centerline {\tir {Complete Conjugacy Invariants of Nonlinearizable
Holomorphic Dynamics}}}

\bigskip
\bigskip

{\centerline {Kingshook Biswas \footnote {$^1$} {Ramakrishna Mission Vivekananda University,
Belur Math, WB-711202, India.
E-mail: kingshook@rkmvu.ac.in} }}

\bigskip
\bigskip

{\bf Abstract.} {\it Perez-Marco proved the existence of
non-trivial totally invariant connected compacts called hedgehogs near the fixed point of a
nonlinearizable germ of holomorphic diffeomorphism. We show that if two nonlinearisable
holomorphic germs with a common indifferent fixed point have a
common hedgehog then they must commute. This allows us to
establish a correspondence between hedgehogs and nonlinearizable
maximal abelian subgroups of Diff$(\bf C,0)$. We also show that two nonlinearizable germs are
conjugate if and only if their rotation numbers are equal and a
hedgehog of one can be mapped conformally onto a hedgehog of the
other. Thus the conjugacy class of a nonlinearizable germ is
completely determined by its rotation number and the conformal
class of its hedgehogs.}

\medskip

{\centerline {\it AMS Subject Classification: 37F50}}

\bigskip

\bigskip

\stit {1. Introduction.}

\bigskip

We consider the dynamics of a holomorphic germ $f(z) = e^{2\pi i
\alpha} z + O(z^2), \, \alpha \in \dd{R} - \dd{Q}$ near the
indifferent irrational fixed point 0. The germ $f$ is said to be {\it
linearizable} if there is a holomorphic change of variables $z =
h(w) = w + O(w^2)$ such that $$ h^{-1} \circ f \circ h = R_\alpha,
$$ where $R_\alpha(w) = e^{2\pi i \alpha} w$ is the rigid
rotation. The maximal linearization domain of $f$ is called the
Siegel disk of $f$.

\medskip

The problem of {\it linearization}, or determining when $f$ is linearizable,
is intimately linked to the arithmetic of the rotation number $\alpha$, and has a long and interesting
history.The work of H. Cremer ([Cr1], [Cr2]) in the 1920's showed the
existence of nonlinearizable germs for rotation numbers very well approximable
by rationals, while that of C.L.Siegel ([Si]) in 1942 and A.D.Brjuno in the
1960's showed linearization was always possible for for germs with rotation
numbers poorly approximated by rationals. The matter was settled definitively
by J.C.Yoccoz ([Yo]) in 1987, when he showed that Brjuno's arithmetic condition
was the optimal one for linearizability. The reader is referred to
R.Perez-Marco's Bourbaki Seminar [PM4] for a complete account of the story.

\medskip

When $f$ is linearisable, the closures of the
linearization domains $h( \{ |w| < r \})$  for $r$ small are
completely invariant connected compacts for $f$. It is not
obvious however whether nonlinearizable germs have completely invariant
non-trivial connected compacts, but Perez-Marco showed
([PM1]) that in fact for any germ $f$ there are always
completely invariant, non-trivial connected compacts near the fixed
point. These invariant compacts $K$ are called {\it Siegel
compacta}. If $K$ is not contained in the closure of a
linearisation domain it is called a {\it hedgehog}. The {\it
hedgehog} is called linearizable or non-linearizable depending on
whether it contains a linearization domain or not.

\medskip

Perez-Marco has studied the topology (in [PM2]) and the dynamics
(in [PM3]) of hedgehogs. The results in [PM2] show that the
topology of hedgehogs is complex. Nonlinearizable hedgehogs have empty interior. They
are not locally connected at any point except possibly at the
fixed point. They always contain points inaccessible from their
complement.

\medskip

The results in [PM3] show that the dynamics on hedgehogs
has many features in common with linearizable dynamics.
For example, there are no periodic points on the hedgehog,
and every point is recurrent (as is the case in the linearizable
situation). Any hedgehog contains a continuous nested
one-parameter family of sub-hedgehogs (recalling the filtration of
a Siegel disk by invariant sub-disks). If two germs commute then
they preserve the same hedgehogs (commuting linearizable germs
have common linearization domains).

\medskip

These results suggest that hedgehogs should be thought of as
"degenerate linearization domains" in some sense, though we are
unable at present to give this heuristic notion a precise
mathematical formulation. Nevertheless it is one of the main
motivations behind Perez-Marco's results and the results of this
article. We prove the converse of the last result mentioned above.
If two germs preserve a common hedgehog then they must commute
(recall that linearizable germs with a common linearization
domain must commute). So the group of germs preserving
the hedgehog is commutative, and equal to the centralizer in
Diff$(\dd C , 0)$ of any element of the group. As Perez-Marco
shows in [PM3], there is a unique continuous nested one-parameter
family of hedgehogs associated to a nonlinearisable germ, indeed
to its centralizer since commuting germs have the same hedgehogs.
On the other hand since germs preserving the same hedgehogs
commute, to every family of hedgehogs we can associate an
abelian subgroup of germs. We show that this gives a one-to-one correspondence between
hedgehogs and the nonlinearizable maximal abelian subgroups of Diff$(\dd C , 0)$
(nonlinearizable meaning not conjugate to the subgroup
$(R_{\lambda}(z) = \lambda z)_{\lambda \in \dd C^*}$).

\medskip

The correspondence is natural with respect to change of variables:
two such subgroups are conjugate by an element $\phi$ of
Diff$(\dd C,0)$ if and only if $\phi$ maps one family of hedgehogs
to the other. Two nonlinearisable germs are conjugate if and only
they have the same rotation number and a hedgehog of one can be
'conformally mapped' to the other (ie there is a conformal mapping
between neighbourhoods of the hedgehogs taking one hedgehog to the
other). Thus the geometry of the hedgehog (in some sense its
conformal geometry) completely determines the nonlinearizable dynamics.
It would be interesting (and morally satisfying) to find an intrinsic notion of
conformal structure for singular spaces like hedgehogs characterizing their
conformal equivalence. It is interesting to compare this with the description of
conjugacy classes of germs tangent to the identity given by
Martinet-Ramis in [Ma-Ra]. They describe a complete set of
invariants as a formal invariant $\lambda$ (a complex number)
together with a singular quotient
object, the "chapelet de spheres", and define an appropriate notion
of conformal equivalence of these objects, so that the conjugacy class
of a germ is determined by the conformal class of the associated chapelet.

\bigskip

\stit {2. Preliminaries.}

\bigskip

Here we collect the basic definitions and results about hedgehogs
due to Perez-Marco which we will require. We let $f(z) = e^{2\pi i
\alpha}z + O(z^2), \alpha \in ${\bf R}, be a germ of holomorphic
diffeomorphism with an indifferent fixed point at $0$.

\medskip

\noindent {\bf Definition 2.1 (Admissible Domain).} An
{\it admissible domain} for $f$  is a Jordan domain with $C^1$-boundary
$U$ containing $0$ such that $f$ and $f^{-1}$ extend univalently
to a neighbourhood of $\overline{U}$.

\medskip

\noindent {\bf Definition 2.2 (Siegel Compacta, Hedgehogs).} A
{\it Siegel compact} of $f$ is a full, connected, compact set
$K$ strictly containing $0$ such that $f$ and $f^{-1}$ extend
univalently to a neighbourhood of $K$ and leave $K$ invariant,
$f(K) = f^{-1}(K) = K$. When the rotation number $\alpha$ is
irrational, a Siegel compact which is not contained in the
closure of a linearization domain of $f$ is called a {\it hedgehog}. A
hedgehog is called linearizable if it contains a linearization
domain and nonlinearizable otherwise.

\medskip

\noindent {\bf Theorem 2.3 (Existence and Uniqueness of Siegel compacta) ([PM1]).}
{\it For any admissible domain $U$ there is a Siegel compact
$K$ contained in $\overline{U}$ which extends upto the boundary
of $U$, i.e. $K \cap \partial U \neq \phi$. When the rotation number $\alpha$
is irrational then there is a unique such Siegel compact $K = K(U)$
which is in fact equal to the connected component
containing $0$ of the set of non-escaping points $\{ z \in \overline{U}: f^n(z)
\in \overline{U} \ \forall n \in {\bf Z} \}$ of $U$. We call $K(U)$
the Siegel compact associated to $U$.}

\medskip

We now restrict ourselves to the case of germs with irrational
rotation number (so any Siegel compact is either a linearization
domain,
or a linearizable hedgehog, or a nonlinearizable hedgehog).
If two admissible domains are nested, $U \subset
V$, then the non-escaping points of $U$ are of course non-escaping
points of $V$, and it follows from the above Theorem that the
associated Siegel compacta are nested, $K(U) \subset K(V)$. In
fact any Siegel compact is filtered by a nested family of
sub-Siegel compacta:


\medskip

\noindent {\bf Theorem 2.4 ([PM3]).}{\it Let $K$ be a Siegel
compact. Given an admissible neighbourhood $U$ such that $K = K(U)$
and a continuous monotone increasing one-parameter family of admissible neighbourhoods $(U_t)_{0
< t \leq 1}, \cap_t U_t = \{ 0 \}, \cup_t U_t = U_1 = U$, the
associated family of Siegel compacta $(K_t = K(U_t))_{0 < t \leq
1}$ is a continuous (for the Hausdorff topology on compact sets),
monotone increasing family of sub-Siegel compacta of $K$, such
that $K_t \to \{0\}$ as $t \to 0$ and $K_t \to K_1 = K$ as $t \to
1$. Moreover any Siegel compact contained in $\overline{U}$
(in particular any sub-Siegel compact of $K$) belongs to the family $(K_t)_{0
< t \leq 1}$.}

\medskip

Thus the sub-Siegel compacta of a given Siegel compact form a
continuous, monotone, one-parameter family (which is a trivial fact
for linearization domains, but quite a remarkable one for hedgehogs,
given their complex topological structures).
The parametrization of the
family is not unique; different choices of the admissible domain
$U$ and the filtration $(U_t)_{0 < t \leq 1}$ lead to different
parametrizations. In the case of a linearizable germ and a linearization
domain, the sub-linearization domains admit a
natural parametrization given by the conformal radius,
which is conformally invariant. An interesting problem is to find
such a conformally invariant parametrization for sub-hedgehogs of
a hedgehog; is there an appropriate notion of
conformal radius for hedgehogs?

\medskip

Any linearizable hedgehog of a germ $f$ must contain the Siegel disk of
$f$; however, it can have no other interior points:

\medskip

\noindent {\bf Theorem 2.5 ([PM3])}{\it The interior of a
linearizable hedgehog is equal to the Siegel disk. The interior of
a nonlinearizable hedgehog is empty.}

\medskip

The family of all Siegel compacta of $f$,
$$
{\cal K} (f) = \{ \ K : \ K \ \hbox{is a Siegel compact of } f \
\},
$$
is a linearly ordered family with respect to inclusion:
given Siegel compacta $K_1, K_2$, they are both sub-Siegel
compacta of their union $K_1 \cup K_2$ (which is
full, since by the maximum principle and Theorem 2.5 the complement cannot have any bounded
components) and hence by Theorem 2.4 one must be contained in the
other. For a nonlinearizable germ $f$ the union $\cup_{K \in {\cal K} (f) } K$
of all its hedgehogs is, morally speaking, the analogue of the Siegel
disk; note that the closure of this union is {\it not} a hedgehog,
since every hedgehog is strictly contained in some
admissible domain and hence in a larger Siegel compact.

\medskip

While the family of all Siegel compacta ${\cal K}(f)$ of a germ $f$ is a natural
object to consider, it is not an invariant of the local dynamics; for any
germ there are of course local changes of variables which do not extend
univalently to neighbourhoods of all the Siegel compacta of the germ.
However, since any change of variables is univalent on a neighbourhood of
every sufficiently small Siegel compact, the situation is easily
remedied by considering the germ of the family of Siegel compacta.

\medskip

\noindent {\bf Definition 2.6 (Germ of Siegel compact, Germ of
Hedgehog).} Consider all families $\cal K$ of compacta in the
plane such that $0 \in K$ for all $K \in {\cal K}$, and such that
for any $\epsilon > 0$ there is a $K \in \cal K$ with diam$(K) < \epsilon$.
We define the following equivalence relation on such familes:


\noindent ${\cal K}_1 \sim {\cal K}_2$ if there is an $\epsilon =
\epsilon({\cal K}_1, {\cal K}_2) > 0$ such that for any compact $K$
with diam$(K) < \epsilon, \ K \in {\cal K}_1 \Leftrightarrow K \in {\cal
K}_2$.

\noindent We call an equivalence class $[\cal K]$ a {\it germ of compact}.

\noindent Given a germ
$f(z) = e^{2 \pi i \alpha}z + O(z^2), \alpha \in {\bf R} - {\bf Q}$,
the {\it germ of Siegel compact} of $f$ is defined to be the
germ of compact $[{\cal K}(f)]$ (where ${\cal K}(f)$ is the family of all
Siegel compacta of $f$). If
$f$ is nonlinearizable, we call $[{\cal K}(f)]$ a {\it germ of
hedgehog}.

\medskip

There is a natural action of germs of homeomorphisms of ${\bf C}$
fixing $0$ on germs of compacta,
$$
(\phi , [{\cal K}] \mapsto \phi([{\cal K}])
$$
where $\phi([{\cal K}])$ is defined as follows: it is possible to
pick a representative $\cal K' \sim \cal K$ such that $\phi$ is
defined in a neighbourhood of all $K \in \cal K'$. Set $\phi([{\cal K}])
= [\{ \ \phi(K) : K \in \cal K' \ \}]$. It is easy to see this
gives a well-defined action.

\medskip

The action restricts to an action of Diff$({\bf C},0)$ on
germs of Siegel compacta:

\medskip

\noindent {\bf Proposition 2.7}.{\it Any germ of diffeomorphism
$\phi \in$Diff$({\bf C},0)$ takes germs of Siegel compacts to
germs of Siegel compacts; indeed  $\phi([{\cal K}(f)]) = [{\cal
K}(\phi \circ f \circ \phi^{-1})]$. Thus germs of Siegel
compacts are holomorphic conjugacy invariants.}

\medskip

\noindent {\bf Proof:} Given $\phi$ in Diff$({\bf C},0)$ and a germ of Siegel compact
$[{\cal K}(f)]$, pick a representative $\cal K \subseteq {\cal K}(f)$
of $[{\cal K}(f)]$ such that $\phi$ is univalent in a
neighbourhood of all $K \in \cal K$. Let $\phi(\cal K)$ be the family
of compacta $\{ \ \phi(K) : K \in \cal K \ \}$. Since $\phi$ takes
invariant sets of $f$ to invariant sets of the conjugate $\phi
\circ f \circ \phi^{-1}$, $\phi({\cal K})$ is a family of Siegel
compacta of $\phi \circ f \circ \phi^{-1}$, which moreover
contains all sufficiently small Siegel compacta of $\phi
\circ f \circ \phi^{-1}$, so $\phi({\cal K}) \sim {\cal
K}(\phi \circ f \circ \phi^{-1})$ and $\phi([{\cal K}(f)]) = [\phi({\cal K})] = [{\cal
K}(\phi \circ f \circ \phi^{-1})]$. $\diamond$

\medskip

We can restate the following result from [PM3] in terms of germs
of hedgehogs:

\medskip

\noindent {\bf Theorem 2.8 ([PM3]).}{\it If two nonlinearisable germs $f,g$
commute then any sufficiently small hedgehog of $f$ is also a
hedgehog of $g$, so they define the same germ of hedgehog $[{\cal K}(f)]
= [{\cal K}(g)]$. }

\bigskip

\stit {3. Results.}

\bigskip

The key result from which the others follow easily is the following:

\medskip

\noindent {\bf Theorem 3.1.} {\it Let $K$ be a hedgehog for a nonlinearizable
germ $f$. If a germ $T(z) = z + O(z^2)$ tangent to the identity is univalent
on a neighbourhood of $K$ and $K$ is either forward or backward invariant
for $T$, then $T$ is equal to the identity.}

\medskip

For notational convenience, given quantities $a, b$ which
are either positive
sequences or functions of $z$ near $0$, we write
$a \preceq b$ if $a \leq C b$ for some constant $C > 0$ for all
$n$ sufficiently large (or all $z$ sufficiently small as the case
may be). We write $a \sim b$ if $a \preceq b$ and $b \preceq a$.

\medskip

\noindent {\bf Proof:} Suppose $T \neq id$, then $T(z) = z + c_d z^{d+1} + O(z^{d+2})$
(for some $d \geq 1, c_d \neq 0$) is a
nondegenerate parabolic germ. Let
$U_1 \subset U_2, U_1 \neq U_2$ be admissible domains for $f$ and $T$ such $K \cap \partial
U_1 \neq \phi, K = K(U_1)$. Let $F \subset \overline{U_2}$ be a Siegel compact
of $T$ such that $F \cap \partial U_2 \neq \phi$. As Perez-Marco shows in [PM1],
$F$ is a invariant Fatou flower of $T$, meaning that int$(F) = P_1 \cup \dots
P_d$ where $P_1, \dots, P_d$ (the 'petals' of the 'flower')
are pairwise disjoint Jordan domains invariant under $T$ whose
boundaries intersect only at the origin, and $T_{|P_i}$ is a parabolic
automorphism of $P_i$ having $0$ as the unique fixed point on the boundary.
Pick a point $z_0 \in \hbox{int}(F) - \overline{U_1}$ and a small
ball $B_0$ around $z_0$ contained in int$(F)- \overline{U_1}$.
Since $B_0 \subset P_i$ for some $i$, $T^n(B_0) \to \{ 0 \}$ as $n
\to \pm \infty$.

\medskip

We may assume wlog (considering $T^{-1}$ if necessary) that $T(K) \subset K$.
For $n \geq 0$ let $z_n = T^{-n}(z_0), B_n = T^{-n}(B_0)$. Then points in $B_n$ escape
from $\overline{U_1}$  under $T^n$, while by hypothesis
points of the hedgehog $K$ remain in $K \subset \overline{U_1}$,
so we will arrive at the desired contradiction if we can
show that $B_n \cap K \neq \phi$ for some $n$. We
need the following estimate on the asymptotic size of the $B_n$'s:

\medskip

\noindent {\bf Lemma 3.2.} \quad {\it $d(z_n, \partial B_n) \succeq |z_n|^{d+1}$}

\medskip

\noindent {\bf Proof:} Taking a covering of a neighbourhood of the
origin by attracting and repelling Fatou petals, we see that the
domains $B_n$ converge to $0$ through a repelling petal $P$; there
is a Fatou coordinate $w = \chi(z)$ defined on $P$ which maps $P$
to the right half-plane $\{ \hbox{Re } w > 0 \}$ and conjugates
$T^{-1}$ to the translation $w \mapsto w + 1$. Moreover,
 $\chi$ has an asymptotic expansion of
the form $\chi(z) = c \log(z) + \sum_{n \geq -d} c_n z^n$ with
$c_{-d} \neq 0$ (see for example [Ec]),
so $|\chi'(z)| \sim |z|^{-(d+1)}$ and $|(\chi^{-1})'(w)| \sim |w|^{-1/d - 1}$.
Fix an $n_0$ such that $B_{n_0} \subset P$. Note that for $n >
n_0$, all the domains $\chi(B_n)$ are translates of
$\chi(B_{n_0})$ and have the same constant diameter, so for large $n$, for
all $w \in \chi(B_n)$, $|w| \sim |\chi(z_n)|$, and hence $|(\chi^{-1})'(w)| \sim |(\chi^{-1})'(\chi(z_n))| = |\chi'(z_n)|^{-1}| \sim |z_n|^{d+1}$.
It follows that for all $z \in B_n$, for $n$ large, putting $w = \chi(z)$,
$|\chi'(z)|^{-1} = |(\chi^{-1})'(w)| \sim |z_n|^{d+1}$.
So taking $z_n' \in \partial B_n$ such that
$|z_n - z_n'| = d(z_n, \partial B_n)$, for $n > n_0$ we have,
for some $z_n'' \in B_n$,
$$\eqalign{
d(z_n, \partial B_n) = |z_n - z_n'| & \geq |\chi'(z_n'')|^{-1}|\chi(z_n) - \chi(z_n')| \cr
                                    & \sim |z_n|^{d+1} d(\chi(z_n), \partial
                                    \chi(B_n)) \cr
                                    & = |z_n|^{d+1} d(\chi(z_{n_0}), \partial
                                    \chi(B_{n_0})) \cr
}$$
and the Lemma follows. $\diamond$

\medskip

The points $z_n$ converge slowly to $0$, in the sense that
$|z_{n+1}| / |z_n| = |T^{-1}(z_n)| / |z_n| = 1 + O(z_n^d) \to 1$ as
$n \to \infty$. To prove the Theorem it suffices to prove the
following Proposition:

\medskip

\noindent {\bf Proposition 3.3} {\it Let $K$ be a hedgehog of a
nonlinearizable germ $f$. Let $(z_n)$ be a sequence converging to
$0$ such that for $n$ large enough $|z_{n+1}| \geq \epsilon |z_n|$
for some $\epsilon > 0$ , and $(B_n)$ a sequence of domains such that $z_n \in B_n$
and $d(z, \partial B_n) \succeq |z_n|^{d+1}$ for some $d \geq 1$.
Then for some subsequence $(B_{n_k})$, $K \cap B_{n_k} \neq \phi$ for all large $k$.}

\medskip

\noindent {\bf Proof:} Note that for any change of variables
$\phi(z) = z + O(z^2)$, the sequences $(\phi(z_n))$ and
$(\phi(B_n))$ satisfy the hypotheses of the Proposition, thus we
may assume wlog, that $f$ is of the form $f(z) = e^{2 \pi i
\alpha}z + O(z^N)$ for some large $N$ which we will choose
appropriately in the course of the proof (any irrationally indifferent germ
can always be analytically conjugated to germs tangent to the
rotation $R_{\alpha}$ upto arbitrarily high orders). We need the following two estimates on how long we can
iterate such a germ close to the origin, and how close the orbits stay to that
of the rotation:

\medskip

\noindent {\bf Lemma 3.4} {\it Given $f(z) = e^{2\pi i \alpha}z + O(z^N)$
for some $N \geq 2$, for all $z$ small enough, at least $M(z)$
iterates of $f$ are defined on $z$ where $M(z) = C |z|^{-N+1}$
for some $C > 0$,
and moreover
$$
|f^k(z)| \leq 2|z| \ , \ k=0,\dots,M(z).
$$
}

\medskip

\noindent {\bf Proof:} There are constants $C_1 > 0$ and $\epsilon_0 > 0$ such
that for $|z| \leq \epsilon_0$ we have
$$
|f(z) - e^{2\pi i \alpha}z| \leq C_1 |z|^N.
$$
So for $|z| \leq \epsilon_0$,
$$
|f(z)| \leq |z| + C_1 |z|^N = \phi(|z|)      \qquad -(1)
$$
where $\phi(t) = t + C_1 t^N$. To estimate $|f^k(z)|$ for small
$z$, we estimate $\phi^k(t)$ for $t = |z|$ close to $0$. It is
convenient to conjugate the mapping $t \mapsto \phi(t)$
to a mapping $s \mapsto \tilde{\phi}(s)$, in terms of the variables
$s = 1/t^{N-1}$ and $\tilde{\phi}(s) = 1/{\phi(t)}^{N-1}$
close to $+\infty$. A calculation gives
$$
\tilde{\phi}(s) = s - (N-1)C_1 + O(1/s) \geq s - C_2
$$
for $s \geq s_0$ sufficiently large, for some constants $s_0, C_2$.
It follows that for $s \geq 2s_0$ and $k \leq s/(2C_2)$, we have
$$
\tilde{\phi}^k(s) \geq s - kC_2
$$
In terms of the variables $\phi(t), t$, this means that for $t
\leq t_0$ sufficiently small and $k \leq 1/(2C_2 t^{N-1})$,
$$
\phi^k(t) \leq t (1 - kC_2 t^{N-1})^{-1 \over N-1}
$$
For $k \leq 1/(2C_2 t^{N-1})$, it is easy to see from the above
that $\phi^k(t) \leq 2t$. So for $|z| \leq $min$(t_0,
\epsilon_0/2)$, we have $\phi^k(|z|) \leq \epsilon_0$ for
$k \leq C |z|^{-N+1}$ where $C = 1/(2C_2)$. Since $\phi$ is
monotone increasing, it follows from (1) by induction
that
$$
|f^k(z)| \leq \phi^k(|z|) \leq 2|z| \leq \epsilon_0 \ , k=0,\dots,
C|z|^{-N+1}
$$
so that at least $M(z) = C |z|^{-N+1}$ iterates of $f$ on $z$ are
defined.
$\diamond$

\medskip

\noindent {\bf Lemma 3.5.} {\it For $|z| \leq \epsilon$, we have
$$
|f^k(z) - R^k_{\alpha}(z)| \leq k C_2 |z|^N \ , \ k = 0,\dots,
M(z)
$$
for some $C_2 > 0$.
}

\medskip

{\bf Proof:} Since $|f^k(z)| \leq 2|z| \leq \epsilon_0$ for $|z| \leq
\epsilon, k=0,\dots,M(z)$, letting $z_k = f^k(z)$, we know
$$
|f(z_k) - R_{\alpha}(z_k)| \leq C_1 |z_k|^N \ , \ k = 0,\dots,
M(z)
$$
so
$$\eqalign{
|f^k(z) - R^k_{\alpha}(z)| & = \left |(f(z_{k-1}) - R_{\alpha}(z_{k-1}))
+ \sum_{j=1}^{k-1}(R^j_{\alpha}(z_{k-j}) - R^{j+1}_{\alpha}(z_{k-j-1})) \right| \cr
& \leq |(f(z_{k-1}) - R_{\alpha}(z_{k-1}))| + \sum_{j=1}^{k-1}|(R^j_{\alpha}(z_{k-j}) -
R^{j+1}_{\alpha}(z_{k-j-1}))| \cr
& = |(f(z_{k-1}) - R_{\alpha}(z_{k-1}))| + \sum_{j=1}^{k-1}|(R^j_{\alpha}(g(z_{k-j-1})) -
R^{j}_{\alpha}(R_{\alpha}(z_{k-j-1})))| \cr
& = \sum_{j=0}^{k-1}|f(z_j) - R_{\alpha}(z_j))| \cr
& \leq \sum_{j=0}^{k-1} C_1 |z_j|^N \cr
& \leq k C_2 |z|^N \cr
}
$$
since $|z_j| \leq 2|z|, j \leq k-1 \leq M(z)$. $\diamond$

\medskip

\noindent {\bf Proof of Proposition 3.3:} For all $n$ sufficiently large, the
circle $\{ |z| = |z_n| \}$ intersects the hedgehog $K$ since
it is a non-trivial connected set containing the origin; let $w_n$
be a point of $K$ on this circle. Let $(p_k/q_k)_{k \geq 0}$ be the continued fraction convergents
of $\alpha$. For all $k$, it follows from the continued
fraction algorithm that any point on the circle $\{ |z| = |z_n| \}$ is at
distance at most $2q^{-1}_k 2\pi  |z_n|$ from the first $q_k$
points $R^m(w_n), m=0,\dots,q_k$ of the orbit of $w_n$ under the rotation $R_{\alpha}$.
If for each $k$ we take $z_{n_k}$ to be be the first element of the sequence $(z_n)$ such that
$|z_{n_k}| < q_k^{-1/(d+1)}$, so $q_k^{-1/(d+1)} \leq |z_{n_k - 1}|$ then by hypothesis
$|z_{n_k}| \geq \epsilon |z_{n_k - 1}|
\geq \epsilon q_k^{-1/(d+1)}$, thus $|z_{n_k}| \sim q_k^{-1/(d+1)}$, and $2q^{-1}_k 2\pi  |z_{n_k}|
 \sim |z_{n_k}|^{d+1} |z_{n_k}| = |z_{n_k}|^{d+2}$. So for all $k$ there is
some $0 \leq m_k \leq q_k$ such $|z_{n_k} -
R^{m_k}(w_{n_k})| \preceq |z_{n_k}|^{d+2}$. By Lemma $3.5$, the
orbit of $w_{n_k}$ under $f$ stays close to the orbit under
$R_{\alpha}$ for at least $M(w_{n_k})$ iterates, and $M(w_{n_k})
\sim |w_{n_k}|^{-(N-1)} = |z_{n_k}|^{-(N-1)} \sim q_k^{(N-1)/(d+1)}$,
so assuming $N \geq d+2$ we have $M(w_{n_k}) \geq q_k$ for all
large $k$, and $|f^{m_k}(w_k) - R^{m_k}(w_k)| \leq q_k C_2 |z_{n_k}|^N \sim |z_{n_k}|^{-(d+1)}
|z_{n_k}|^N = |z_{n_k}|^{N - (d+1)}$. Thus by taking $N \geq 2d + 3$ we have
$|f^{m_k}(w_k) - R^{m_k}(w_k)| \preceq |z_{n_k}|^{d+2}$, and
$$\eqalign{
|z_{n_k} - f^{m_k}(w_k)| & \leq |z_{n_k} -
R^{m_k}(w_{n_k})| + |f^{m_k}(w_k) - R^{m_k}(w_k)| \cr
                         & \preceq |z_{n_k}|^{d+2} \cr
}$$
Since $d(z_{n_k}, \partial B_{n_k}) \succeq |z_{n_k}|^{d+1}$, it
follows that $f^{m_k}(w_k) \in B_{n_k}$ for all large $k$ and
$K \cap B_{n_k} \neq \phi$.
$\diamond$

\medskip

\noindent {\bf Proof of Theorem 3.1.} By Proposition 3.3 we can pick a point $w$ in
$K \cap B_n$ for some $n > 0$; then $T^n(w) \in B_0 \subset \overline{U_1}^c$ whereas
by hypothesis $T^n(w) \in K \subset \overline{U_1}$, a contradiction. $\diamond$

\medskip

%
%
%

\noindent {\bf Theorem 3.6.} {\it A nonlinearizable germ $f$ cannot
commute with a nondegenerate parabolic germ $g$.}

\medskip

\noindent {\bf Proof:} While this is easy to see using formal power series
arguments, we can give a dynamical proof using the previous
result. Indeed if $f$ and $g$ commute, then for a small hedgehog $K$ of
$f$ on a neighbourhood of which $g, g^{-1}$ are univalent,
$f(g(K)) = g(f(K)) = g(K)$, so $g(K)$ is also a hedgehog for $f$.
Since the hedgehogs of $f$ are linearly ordered with respect to inclusion,
$K$ is either forward or backward invariant under $g$, contradicting Theorem 3.1,
since $g^q$ is tangent to the identity for some $q$ but $g^q \neq id$. $\diamond$

\medskip

\noindent {\bf Theorem 3.7.} {\it Two nonlinearisable germs $f_1$ and $f_2$
are holomorphically conjugate by a germ $\phi$ if and only if $\phi$
maps some hedgehog $K_1$ of $f_1$ to a hedgehog $K_2$ of $f_2$ and
the rotation numbers of $f_1$ and $f_2$ are equal. Thus for a nonlinearisable
germ $f$ its rotation number and germ of hedgehog are a complete set of
holomorphic conjugacy invariants.}

\medskip

\noindent {\bf Proof:} If $f_2 = \phi \circ f_1 \circ \phi^{-1}$, then for a
hedgehog $K_1$ of $f_1$,
$$
f_2(\phi(K_1)) = \phi \circ f_1 \circ \phi^{-1} (\phi(K_1)) =
\phi(f_1(K_1)) = \phi(K_1)
$$
which implies that $\phi(K_1)$ is a hedgehog of $f_2$. This proves
the "only if" part.

\medskip

For the "if" part, given $\phi(K_1) = K_2$, let $\tilde{f}_2$
be the conjugate $\tilde{f}_2 := \phi \circ f_1 \circ \phi^{-1}$.
Then $\tilde{f}_2(K_2) = K_2$, so the germ $g = \tilde{f}_2 \circ
f^{-1}_2$ preserves $K_2$. Since the rotation numbers of $f_1$ and
$f_2$ are equal, $g$ is tangent to the identity and hence by Theorem 3.1
is equal to the identity. So
$\tilde{f}_2 = f_2$, i.e. $\phi$ conjugates $f_1$ to $f_2$. $\diamond$

\medskip

We have the converse of Theorem 2.8:

\medskip

\noindent {\bf Theorem 3.8.} {\it If two nonlinearizable germs $f$ and $g$
have a common hedgehog $K$ then they commute. In particular if $f$ and $g$
have the same germ of hedgehog (so all sufficiently small hedgehogs of one are
hedgehogs of the other) then they commute.}

\medskip

\noindent{\bf Proof:} The hedgehog $K$ is invariant under the
commutator $f \circ g \circ f^{-1} \circ g^{-1}$ which is
tangent to the identity and therefore by Theorem 3.1 must be the
identity. $\diamond$

\medskip

We now consider abelian subgroups $H$ of Diff$({\bf C},0)$
such that $H$ contains at least one irrationally indifferent germ. If any
irrationally indifferent germ in $H$ is linearizable then all
germs in $H$ are linearizable and we call $H$ linearizable,
otherwise all irrationally indifferent germs in $H$ are
nonlinearizable (the rationally indifferent germs in $H$ are
of finite order and linearizable) and we call $H$ nonlinearizable.
Given a nonlinearizable abelian
subgroup $H$, by Theorem 2.8 all irrationally indifferent germs in
$H$ determine the {\it same} germ of hedgehog, which we denote by
$[{\cal K}(H)]$. Conversely, to every germ of hedgehog $[\cal K]$, we
can associate the subgroup of germs which leave it invariant,
Aut$([{\cal K}]) := \{ f \in \hbox{Diff}(\dd C,
0) : f([{\cl K}]) = [{\cl K}] \}$.

\medskip

\noindent {\bf Theorem 3.9} {\it For any germ of hedgehog $[\cal
K]$, the homomorphism
$$\eqalign{
\lambda : Aut([{\cal K}]) & \to {\bf C}^* \cr
 g              & \mapsto g'(0) \cr
 }$$
 is an injective homomorphism into $S^1 \subset {\bf C}^*$. Thus
 Aut$([{\cal K}])$ is abelian and germs in Aut$([{\cal K}])$ are
 uniquely determined by their rotation numbers.}

\medskip

\noindent{\bf Proof:} For any $g \in $Aut$([\cal K])$, for all sufficiently
small $K \in \cal K$, $g(K) \in \cal K$, and since hedgehogs of a nonlinearizable
germ are linearly ordered, $K$ is either forward or backward invariant
under $g$. Any $g$ in the kernel of $\lambda$
is tangent to the identity and by Theorem 3.1 equals the
identity. Thus $\lambda$ is injective and Aut$([\cal K])$ is abelian.
If $|g'(0)| \neq 1$ for some $g \in $Aut$([\cal K])$ then $g$ is linearizable,
and since Aut$([\cal K])$ is abelian, the linearization of $g$ linearizes
all germs in Aut$([\cal K])$, contradicting the existence of at least one
nonlinearizable germ in Aut$([\cal K])$. Thus $|g'(0)| = 1$ for all
$g \in $Aut$([\cal K])$. $\diamond$

\medskip

Finally we have:

\medskip

\noindent {\bf Theorem 3.10.} {\it There is a bijective correspondence
between nonlinearizable maximal abelian
subgroups of Diff$(\dd C, 0)$ and germs of hedgehogs:
$$\eqalign{
H & \rightarrow [{\cal K}(H)] \cr
Aut([\cl K])  & \leftarrow [\cal K] \cr
}$$
The action of Diff$(\dd C, 0)$ on nonlinearisable
maximal abelian subgroups by conjugation $(h, H) \mapsto h \, H \,
h^{-1}$ corresponds to the action of Diff$(\dd C, 0)$ on
germs of hedgehogs $(h , [\cl K]) \mapsto h([\cl K])$
: $H_1 = h \, H_2 \, h^{-1}$, if and only if $h([\cl K(H_1)]) = [\cl K(H_2)]$.
Thus the germ of hedgehog $[{\cal K}(H)]$ of a nonlinearizable maximal
abelian subgroup $H$ is a complete conjugacy
invariant of $H$.}

\medskip

\noindent {\bf Proof:} We first check that every subgroup Aut$([\cl K])$ is
maximal abelian. Let $f \in $Aut$([{\cal K}])$ be a nonlinearizable germ such
that $[{\cal K}(f)] = [\cal K]$. Any germ $g$ commuting with Aut$([\cl K])$
commutes with $f$ and hence maps sufficiently small
hedgehogs $K$ of $f$ to hedgehogs of $f$ since $f(g(K)) = g(f(K)) = g(K)$, thus
$g([\cal K]) = [\cal K]$ and $g \in $Aut$([{\cal K}])$.


\medskip

It is straightforward to check that the two maps are mutual
inverses. The second assertion of the Theorem follows easily from
Theorem 3.7. $\diamond$

\medskip

\noindent {\bf Remark.} It is well known that the unique maximal
abelian subgroup containing an irrationally indifferent germ $f$
is its centralizer Cent$(f)$, thus Aut$([\cal K]) = $Cent$(f)$
for any nonlinearizable $f$ whose germ of hedgehog is $[\cal K]$.

\bigskip

\noindent {\bf 5. Acknowledgements.} I thank Ricardo Perez-Marco for
suggesting the method of proof of Theorem 3.8 (via consideration of
the parabolic commutator germ and Theorem 3.1), and for numerous
helpful discussions about hedgehogs. I am
also grateful to Arnaud Cheritat for suggesting the proof of
Proposition 3.3 via coordinate changes and elementary estimates;
the original proof was based on Yoccoz's renormalization for germs,
which gives much sharper but non-elementary estimates. I thank Xavier Buff as
well for many helpful discussions.

\bigskip

\centerline{\bf References}

\bigskip

\noindent [Cr1] H.CREMER, {\it Zum Zentrumproblem}, Math Ann., vol 98, 1928, p. 151-153

\medskip

\noindent [Cr2] H.CREMER, {\it \"Uber die H\"aufigkeit der Nichtzentren}, Math Ann., vol 115, 1938, p.573-580

\medskip

\noindent [Ec] J.ECALLE, {\it Th\'eorie it\'erative: introduction
\`a la th\'eorie des invariantes holomorphes}, J. Math pures et
appliqu\'es {\bf 54}, 1975, p. 183-258.

\medskip

\noindent [Ma-Ra] J.MARTINET \& J.P.RAMIS, {\it Classification analytique des
\'equations differentielles non lin\'eaires r\'esonnantes du
premier ordre}, Ann. Sci. \'Ecole Norm. Sup. (4) {\bf 16} (1983),
no. 4, p. 571-621 (1984).

\medskip

%

\noindent [PM1] R.PEREZ-MARCO, {\it Fixed points and circle maps}, Acta Mathematica, {\bf 179:2}, 1997, p.243-294.

\medskip

\noindent [PM2] R.PEREZ-MARCO, {\it Topology of Julia sets and hedgehogs}, Preprint, Universit\'e de PARIS-SUD, 1994.

\medskip

\noindent [PM3] R.PEREZ-MARCO, {\it Hedgehog's Dynamics}, Preprint, UCLA, 1996.

\medskip

\noindent [PM4] R.PEREZ-MARCO, {\it Solution compl\`ete au
probl\`eme de Siegel de lin\'earisation d'une application
holomorphe au voisinage d'un point fixe (d'apr\`es J.C.Yoccoz)},
Seminar Bourbaki {\bf 753}, 1991.

\medskip

\noindent [Si] C.L.SIEGEL, {\it Iteration of Analytic Functions}, Ann. Math., vol.43, 1942, p. 807-812

\medskip

\noindent [Yo] J.C.YOCCOZ, {\it Petits diviseurs en dimension 1}, S.M.F., Asterisque {\bf 231} (1995)

\medskip
\medskip

\noindent Ramakrishna Mission Vivekananda University,
Belur Math, WB-711202, India.

\noindent Email: kingshook@rkmvu.ac.in

\end